\documentclass[11pt]{amsart}
\usepackage{amscd,amssymb,verbatim}
\theoremstyle{plain}
\newtheorem{Thm}{Theorem}[section]
\newtheorem{Cor}[Thm]{Corollary}
\newtheorem{Lem}[Thm]{Lemma}
\newtheorem{Prop}[Thm]{Proposition}
\newtheorem{Def}[Thm]{Definition}
\newtheorem{remark}{Remark}

\newtheorem{notation}{Notation}

\begin{document}
\title[A continuum of totally incomparable H.I. spaces]
{A continuum of totally incomparable Hereditarily 
indecomposable Banach spaces}
\author{I. Gasparis}
\address{Department of Mathematics \\
         Oklahoma State University \\
         Stillwater, OK 74078-1058 \\
         U.S.A.}
\email{ioagaspa@math.okstate.edu}
\begin{abstract}
A family is constructed of cardinality equal to the continuum,
whose members are totally incomparable Hereditarily Indecomposable
Banach spaces.
\end{abstract}
\maketitle
\section{Introduction} \label{S:1}
All Banach spaces considered in this paper are real,
infinite dimensional. By a {\em subspace} of a Banach space
we shall mean an infinite dimensional, closed linear subspace. 
A Banach space is said to be 
{\em Hereditarily Indecomposable} (H.I.) if for every pair
\(Y\), \(Z\) of subspaces
of \(X\) with \(Y \cap Z =\{\mathbf{0}\}\), the subspace
\(Y + Z\) is not closed.
The famous example of Gowers and Maurey \cite{gm} of a Banach
space without unconditional basic sequence, was observed by
W. Johnson to be H.I. Since the appearance of the Gowers-Maurey
space the study of H.I. spaces has been one of the most important
research topics in modern Banach space theory. We refer to
\cite{o} and \cite{af} for a detailed survey of results.

It is proved in \cite{af} that every Banach space not containing
an isomorph of \(\ell_1\) has a subspace which is a quotient
of an H.I. space. A recent result of S. Argyros \cite{ar} states
that a separable Banach space universal for the class of
reflexive H.I. spaces, is also universal for the class of 
separable Banach spaces.
Both results indicate the large variety of H.I.
spaces. The aim of this paper is towards this direction.
Our main result is the following:
\begin{Thm} \label{T:1}
There exists a family of cardinality equal to the continuum
whose members are totally incomparable, reflexive H.I. spaces.
\end{Thm}
Recall that the Banach spaces \(X\) and \(Y\) are {\em totally
incomparable} if no subspace of \(X\) is isomorphic to a subspace
of \(Y\).

The construction of H.I. spaces is not an easy task.
The crucial step was Schlumprecht's construction of an 
arbitrarily distortable Banach space \cite{sch}.
Recall that the Banach space \((X, \| \cdot \|)\)
is {\em arbitrarily distortable} if for every \(\lambda > 1\),
there exists an equivalent norm \(|\cdot|\) on \(X\) 
so that for every subspace
\(Y\) of \(X\) there exist non-zero vectors \(x\), \(y\) in \(Y\)
such that \(\|x\|=\|y\|\), yet \(|x|/|y| > \lambda\).
Schlumprecht's space had an immense
impact in the development of the theory because of its connection
to the Gowers-Maurey construction, as well as to the solution 
of the distortion
problem for \(\ell_p\), \(1 < p < \infty\), \cite{os}.

The first example of an arbitrarily distortable, asymptotic
\(\ell_1\) space was given in \cite{ad}. They showed that
there exist infinite subsets \(M=(m_i)\), \(N=(n_i)\) of 
\(\mathbb{N}\) so that the mixed
Tsirelson space
\(T(\frac{1}{m_i}, S_{n_i})_{i=1}^{\infty}\), is arbitrarily
distortable. In the same paper this example was conditionalized 
to yield an asymptotic \(\ell_1\) H.I. space.

The proof of Theorem \ref{T:1} is based on ideas from \cite{ad}.
However, our argument is considerably simpler. We shall next describe
how this paper is organized. In Section \ref{S:3} we introduce,
for a given scalar \(d > 1\), the infinite subsets \(N\) and
\(P\) of \(\mathbb{N}\) and the null scalar sequence \(\mathbf{a}\),
the \((d,N,P,\mathbf{a})\) distortion property, Definition 
\ref{dp}, of an asymptotic \(\ell_1\) Banach space. This property
will enable us to give a criterion, Theorem \ref{tdp},
for an an asymptotic \(\ell_1\) Banach space to be arbitrarily
distortable. We also show how to obtain totally incomparable
arbitrarily distortable spaces. We apply Theorem \ref{tdp}
in Section \ref{S:4} in order to give an alternative proof
of the fact that certain mixed Tsirelson spaces are
arbitrarily distortable \cite{ad}, \cite{ano}, \cite{adm}.
These spaces can be described as the completion of
\(c_{00}\), the space of all ultimately vanishing real sequences,
under the norm given by \(\|x\|= \sup \{
\sum_{i=1}^{\infty} \mu (\{i\}) x(i) : \mu \in \mathcal{M} \}\),
where \(\mathcal{M}\) is a suitable symmetric subset of the 
finitely supported signed measures on \(\mathbb{N}\)
containing the point mass measures and closed under interval
restrictions. The main difficulty in the study of mixed Tsirelson
spaces is that the norming set \(\mathcal{M}\) is defined by means
of an inductive procedure. We are able to by pass this difficulty
by describing \(\mathcal{M}\) analytically and proving a decomposition
result for its members, Lemma \ref{L:2},  
which greatly simplifies the argument for the distortion of  
\(T(\frac{1}{m_i}, S_{n_i})_{i=1}^{\infty}\). 

In Section \ref{S:5}, we choose a subset
\(\mathcal{N}\) of \(\mathcal{M}\) which is maximal 
with respect to a Maurey-Rosenthal type of condition
\cite{mr} and show in Theorem \ref{dhi} that the completion of 
\(c_{00}\) under the norm induced by \(\mathcal{N}\) is an H.I.
space satisfying a \((d,N,P,\mathbf{a})\) distortion property.
Various choices of \(\mathcal{N}\) give rise to totally incomparable
H.I. spaces.

In order to prove that a space \(X\) is H.I., we employ Theorem
\ref{hi} which loosely speaking asserts that if for every
\(\epsilon > 0\) there exist integers \(k < n\) such that 
every block subspace \(Y\) of \(X\) contains a sufficiently large
(in the Schreier sense) block basis
\(z_1 < \cdots < z_p\) with the property that
\(\|\sum_{i=1}^p a_i z_i\| \geq \epsilon
  \|\sum_{i=1}^p a_i e_i\|_n\),
whenever \((a_i)_{i=1}^p \subset \mathbb{R}^+\),
while 
\(\|\sum_{i=1}^p a_i z_i\| \leq 
  \|\sum_{i=1}^p a_i e_i\|_{Ck}\),   
for every sequence \((a_i)_{i=1}^p\) in \(\mathbb{R}\),
then \(X\) contains no infinite unconditional sequence.
In the above, \((e_i)\) is the natural unit vector basis of
\(c_{00}\) and \(\| \cdot \|_n\), \(\| \cdot \|_{Ck}\)
denote the \(n\)-th Schreier and \(k\)-th conditional Schreier
norms respectively.

The precise statements for the results mentioned above are given
in Section \ref{S:3}. The proof of Theorem \ref{T:1}, presented in
Section \ref{S:3}, follows from Theorem \ref{dhi} and Proposition
\ref{tin}
combined with two fundamental results of descriptive set theory,
the infinite Ramsey theorem \cite{el}, \cite{o2}
and a theorem of Kuratowski \cite{ku}.
\section{Preliminaries} \label{S:2} 
We shall make use of standard Banach space facts and terminology
as may be found in \cite{lt}. If \(D\) is any set, we let
\([D]\) (resp. \([D]^{< \infty}\)) denote the set of its infinite
(resp. finite) subsets. 
Given \(M \in [\mathbb{N}]\), the notation \(M=(m_i)\)
indicates that \(M=\{m_1 < m_2 < \cdots \}\).
Let \(E\) and \(F\) be finite subsets of \(\mathbb{N}\).
We write \(E < F\) if \(\max E < \min F\).  

Suppose now that \(X\) is a Banach space with a Schauder basis
\((e_n)\). A sequence \((u_n)\) in \(X\) is a {\em block basis}
of \((e_n)\) if there exist successive subsets \(F_1 < F_2 < \cdots\)
of \(\mathbb{N}\) and a scalar sequence \((a_n)\) so that
\(u_n = \sum_{ i \in F_n} a_i e_i\), for every
\(n \in \mathbb{N}\). We adopt the notation
\(u_1 < u_2 < \cdots \) to indicate that \((u_n)\)
is a block basis of \((e_n)\). We let
\(\mathrm{supp} u_n\) denote the set \(\{i \in F_n: a_i \ne 0 \}\).  
The {\em range} \(r(u_n)\) of \(u_n\), is the smallest
integer interval containing \(\mathrm{supp} u_n\). The subspace of
\(X\) generated by a block basis of \((e_n)\) is called a 
{\em block subspace}.
We next review two important hierarchies. The Schreier
hierarchy
\(\{S_\xi\}_{\xi < \omega_1}\),
\cite{aa} and the repeated averages hierarchy,
\((\xi_n^M)_{n=1}^\infty\), \(\xi < \omega_1\), \(M \in [\mathbb{N}]\),
\cite{amt}.
Since we shall only be using the families
\(\{S_\xi\}_{\xi < \omega}\), and  
\((\xi_n^M)_{n=1}^\infty\), \(\xi < \omega\), 
\(M \in [\mathbb{N}]\), we confine the definitions to the finite
ordinal case. 

{\bf The Schreier families}. We let 
\(S_0 = \bigl \{ \{n\} : n \in \mathbb{N} \} \bigr \} \cup
               \{\emptyset\}\).
Suppose \(S_\xi\) has been defined, \(\xi < \omega\).
We set 
\[ 
S_{\xi + 1} =   \{ \cup_{i=1}^n F_i : \, 
                     n \in \mathbb{N}, \, n \leq \min F_1,  
             F_1 < \cdots < F_n, \,      
                     F_i \in S_\xi \, (i \leq n) \}
        \cup \{\emptyset\}. 
\]
An important property shared by the Schreier families
is that they are {\em hereditary}: If \(F \in S_\xi\)
and \(G \subset F\), then \(G \in S_\xi\). Another important property
is that they are spreading: If \(\{p_1, \cdots , p_k\}
\in S_\xi\), \(p_1 < \cdots < p_k\), and \(q_1 < \cdots < q_k\)
are so that \(p_i \leq q_i\) for all \(i \leq k\), then
\(\{q_1, \cdots , q_k\} \in S_\xi\).
It is not hard to check that if \(F_1 < \cdots < F_n\)
are members of \(S_\alpha\) such that 
\(\{\min F_i : i \leq n\}\) belongs to \(S_\beta\), then
\(\cup_{i=1}^n F_i\) belongs to \(S_{\alpha + \beta}\).

{\bf The repeated averages hierarchy}.
We first let \((e_n)\) denote the unit vector
basis of \(c_{00}\).  
Given \(\xi < \omega\) and 
\(M \in [\mathbb{N}]\), we define by induction, a sequence 
\(({\xi}_{n}^{M})_{n=1}^{\infty}\) of finitely supported probability
measures on \(\mathbb{N}\) whose supports are successive subsets
of \(M\).

If \(\xi = 0\), then \({\xi}_{n}^{M} = e_{m_n}\), for all 
\(n \in \mathbb{N}\), where \(M = (m_n)\). 

Assume that \(({\xi}_{n}^{M})_{n=1}^{\infty}\)
has been defined for all \(M \in [\mathbb{N}]\).
Set
\[ {[\xi + 1]}_{1}^{M} = \frac{1}{m_1}
                   \sum_{i=1}^{m_1} {\xi}_{i}^{M}
\]
where \(m_1 = \min M\). Suppose that
\( {[\xi + 1]}_{1}^{M} < \cdots < {[\xi + 1]}_{n}^{M}\)
have been defined. Let
\[
   M_n = \{ m \in M : m > \max \, \mathrm{supp}\, {[\xi + 1]}_{n}^{M} \}
\text{ and } k_n = \min M_n .\]
Set
\[{[\xi + 1]}_{n + 1}^{M} = \frac{1}{k_n} \sum_{i=1}^{k_n}
  {\xi}_{i}^{M_n} .
\]
It follows that \(\mathrm{supp} \, \xi_n^M\) belongs to
\(S_\xi\), and moreover it is a maximal (under inclusion) member of
\(S_\xi\). It can be easily shown, by induction, that
if \(i\) and \(j\) belong to \(\mathrm{supp} \, \xi_n^M\)
and \(i < j\), then \(\xi_n^M (\{i\})
\geq \xi_n^M (\{j\})\).

For a probability measure \(\mu\) in \(\mathbb{N}\) and
\(\xi < \omega\), we set \(\|\mu\|_\xi =
\sup \{ \mu (F) : F \in S_\xi\}\). It is proven in \cite{gl},
\cite{ag} that \(\|\xi_1^M\|_{\xi -1} \leq \frac{\xi}{\min M}\),
for every \(\xi \geq 1\) and \(M \in [\mathbb{N}]\).
It follows that for every \(P \in [\mathbb{N}]\), every \(\xi \geq 1\) 
and every \(\epsilon > 0\), there exists \(M \in [P]\)
such that \(\|\xi_1^M\|_{\xi -1} < \epsilon\).
This property of the repeated averages will be very useful in the sequel.
For a detailed study of these hierarchies we refer to
\cite{aa}, \cite{amt}, \cite{otw}, \cite{g}, \cite{ag} and \cite{gl}.

We continue by introducing some more terminology.
A finite collection \(\mathcal{F}\)
of finite subsets of \(\mathbb{N}\) is said
to be \(r S_\xi\)-admissible, \(\xi < \omega\), 
\(r \in \mathbb{N}\), if
there exists an enumeration \(\{I_k : k \leq n\}\)
of \(\mathcal{F}\) such that \(I_1 < \cdots < I_n\)
and the set \(\{\min I_k : k \leq n \}\) 
is the union of \(r\) members of \(S_\xi\).
In case \(\{\min I_k : k \leq n \}\)
is a maximal (under inclusion) member of \(S_\xi\),
\(\mathcal{F}\) is called maximally \(S_\xi\)-admissible.
A finite block basis \(u_1 < \cdots < u_n\) in a Banach
space with a basis is \(r S_\xi\) (resp.
maximally \(S_\xi\))-admissible, if 
\(\{\mathrm{supp} \, u_i : i \leq n\}\) is.

In what follows, \(X\) is a Banach space with a basis \((e_n)\).
The support of every block basis of \((e_n)\) will always
be taken with respect to \((e_n)\).
\begin{Def}
Let \((u_n)\) be a normalized block basis of \((e_n)\), \(\epsilon > 0\)
and \(1 \leq \xi < \omega\). Set \(p_n= \min \mathrm{supp} \, u_n\),
\(n \in \mathbb{N}\), and \(P= (p_n)\). 
\begin{enumerate}
\item A generic \((\epsilon, \xi)\)
average of \((u_n)\) is any vector of the form 
\(\sum_{n=1}^\infty \xi_1^R (p_n) u_n\), where
\(R \in [P]\) and \(\|\xi_1^R\|_{\xi - 1} < \epsilon\).
\item An \((\epsilon, \xi)\) average of \((u_n)\) is any generic
\((\epsilon, \xi)\) average of a normalized block basis
of \((u_n)\). 
\item A normalized \((\epsilon, \xi)\) average of \((u_n)\)
is any vector \(u\) of the form
\(u= \frac{v}{\|v\|}\), where \(v\) is an 
\((\epsilon, \xi)\) average of \((u_n)\). In case \(\|v\| \geq 
\frac{1}{2}\), \(u\) is a smoothly normalized
\((\epsilon, \xi)\) average of \((u_n)\).
\end{enumerate}
\end{Def}
\begin{notation}
Let \(E^*\) be a finite collection of successive intervals of
\(\mathbb{N}\) and let \(u\) be a finite linear combination of
\((e_n)\). 
\begin{enumerate}
\item We let \(I(u,E^*)\) denote the number of elements of
\(E^*\) which are intersected by \(\mathrm{supp} \, u\).
\item Let \(D\) be a finite block basis of \((e_n)\)
such that the support of every member of \(D\) intersects
at least one member of \(E^*\). We set
\(D(E^* , 1)=\{ u \in D: I(u,E^*)=1\}\) and 
\(D(E^* , 2)=\{ u \in D: I(u,E^*) \geq 2\}\).
\end{enumerate}
\end{notation}
Before closing this section, we recall the definitions of
the Schreier space, \(X^\xi\), and conditional Schreier
space, \(CX^\xi\), \(\xi < \omega\).
\(X^\xi\) is the completion of \(c_{00}\) under the norm 
\(\|x\|_\xi= \sup \{\sum_{i \in F} |x(i)| : F \in S_\xi\}\).
\(X^0\) is isometric to \(c_0\).
\(X^1\) was introduced by Schreier \cite{sc} in order to
provide an example of a weakly null sequence without
Cesaro summable subsequence. The generalized family of
Schreier spaces \(\{X^\xi\}_{\xi < \omega_1}\)
was studied in \cite{aa}, where it is shown that
the natural Schauder basis \((e_n)\) of \(X^\xi\)
is \(1\)-unconditional and shrinking.
For a detailed study of the spaces \(\{X^\xi\}_{\xi < \omega}\)
we refer to \cite{gl}.

The conditional Schreier spaces \(\{CX^\xi\}_{\xi < \omega}\),
were constructed by H. Rosenthal (unpublished).
\(CX^\xi\) is the completion of \(c_{00}\) under
the norm 
\[\|x\|_{C\xi}= \sup \biggl \{ \sum_{k=1}^n \biggl |
\sum_{i \in J_k} x(i) \biggr | : n \in \mathbb{N}, \, (J_k)_{k=1}^n 
\text{ are } S_\xi \text{ admissible intervals }
\biggr \}.\]
The natural basis \((e_n)\) of \(CX^\xi\) is of course,
a conditional basis. When \(\xi=0\), \((e_n)\) is equivalent
to the summing basis of \(c_0\).
We also mention the following useful fact:
Suppose \((a_i)_{i=1}^n\) is a non-increasing finite
sequence of non-negative scalars. Then
\(\|\sum_{i=1}^n (-1)^i a_i e_{t_i} \|_{C\xi} \leq
  \|\sum_{i=1}^n a_i e_{t_i} \|_\xi\),
for every increasing sequence of integers \((t_i)_{i=1}^n\).
\section{Main results} \label{S:3}
We start this section by recalling that a normalized sequence
\((x_n)\) in a Banach space is an \(\epsilon\)-\(\ell_1^\xi\)
{\em spreading model}, \(\epsilon > 0\), 
if \(\|\sum_{i \in F} a_i x_i\| \geq \epsilon \sum_{i \in F} 
     |a_i|\), for every \(F \in S_\xi\) and all choices of
scalars \((a_i)_{i \in F}\).

A Banach space \(X\) with a basis \((e_n)\) is {\em asymptotic}
\(\epsilon\)-\(\ell_1^\xi\), \(1 \leq \xi < \omega\),
if every normalized block basis
of \((e_n)\) is an \(\epsilon\)-\(\ell_1^\xi\) spreading model.
\(X\) is asymptotic \(\ell_1\), if it is asymptotic
\(\epsilon\)-\(\ell_1^1\), for some \(\epsilon > 0\) \cite{mmt}.
For an asymptotic \(\epsilon\)-\(\ell_1^\xi\) space \(X\) with
a basis \((e_n)\)
and \(\delta > 0\), we define
\begin{align}
\tau(X, \delta)= \sup &\{ \zeta < \omega: \text{ every normalized block
basis of } (e_n) \notag \\
&\text{ has a subsequence which is a }
\delta - \ell_1^\zeta \text{ spreading model }\}. \notag
\end{align}
Evidently, \(\tau(X, \epsilon) \geq \xi\).
The modulus \(\tau(X, \delta)\) is implicitly defined
in \cite{otw} and \cite{ano}. Of course \(\tau(X, \delta)\)
depends on the choice of the basis \((e_n)\), but it will be clear
from the context which basis is used.
In case \(U\) is a block subspace of \(X\), \(\tau(U, \delta)\)
will be calculated with respect to the block basis that
generates \(U\).
\begin{Def} \label{dp}
Let \(X\) be a Banach space with a basis \((e_i)\).
Let \(N=(n_i)\) and \(P=(p_i)\) be infinite subsets of
\(\mathbb{N}\) such that \(n_{i-1} \leq p_i < \frac{n_i}{2}\),
for every \(i \in \mathbb{N}\). Let \(\mathbf{a}=(\delta_i)\)
be a decreasing null sequence of scalars, and let \(d > 1\).
\(X\) is said to satisfy the \((d,N,P,\mathbf{a})\)
distortion property if for every \(j \in \mathbb{N}\),
\(X\) is an asymptotic
\(\delta_j\)-\(\ell_1^{n_j}\) space such that
\(\tau(U, d \delta_j) < p_j\), for every block subspace
\(U\) of \(X\).
\end{Def}
\begin{Thm} \label{tdp}
Let \((X,\|\cdot\|)\) be a Banach space with a normalized, shrinking,
bimonotone basis \((e_i)\). Suppose that there exist
\(N\), \(P\) in \([\mathbb{N}]\), a scalar sequence
\(\mathbf{a}=(\delta_i)\)
and \(d > 1\) so that \(X\) satisfies the
\((d,N,P,\mathbf{a})\) distortion property.
Then \(X\) is arbitrarily distortable.
\end{Thm}
\begin{proof}
In the sequel the admissibility of every block basis
of \((e_i)\) will always be considered with respect to \((e_i)\).
Given \(j \in \mathbb{N}\), we set
\[\mathcal{A}_j = \biggl \{ \delta_j \sum_{i=1}^k
x_i^* : \, (x_i^*)_{i=1}^k \subset B_{X^*} \text{ is }
    S_{n_j}- \text{admissible} \biggr \}.\]
In the above, the admissibility of \((x_i^*)\) is measured
with respect to \((e_i^*)\), the sequence of functionals
biorthogonal to \((e_i)\). Because \(\tau(X,\delta_j) \geq n_j\),
we have that \(\mathcal{A}_j \subset B_{X^*}\).
Indeed, suppose that \(\delta_j \sum_{i=1}^k x_i^* \in \mathcal{A}_j\)
and let \(x \in X\), \(\|x\| \leq 1\). Put
\(x_i = x | \mathrm{ran}(x_i^*)\), \(i \leq k\).
Since \((e_i)\) is bimonotone, \(\|\sum_{i=1}^k x_i\| \leq 1\).
Furthermore, \((x_i)_{i=1}^k\) is \(S_{n_j}\) admissible. Hence,
\(\delta_j \sum_{i=1}^k \|x_i\| \leq 1\) and the assertion follows.

We define an equivalent norm \(\| \cdot \|_j\) on \(X\) in the following
manner:
\[ \|x\|_j = \delta_j \|x\| + 
\sup \{ x^*(x): \, x^* \in \mathcal{A}_j\}.\]
Let \((u_i)\) be a normalized block basis of \((e_i)\), and let
\(j_0 \in \mathbb{N}\). Let \(U\) be the block subspace of
\(X\) generated by \((u_i)\). Since \(\tau(U, d \delta_{j_0})
< p_{j_0}\),  there exists a normalized block basis 
\((v_i)\) of \((e_i)\) in \(U\) having no subsequence which is
a \(d \delta_{j_0}\)-\(\ell_1^{p_{j_0}}\) spreading model.
It follows, by the main result of \cite{g} combined with
Corollary 3.6 of \cite{ag}, that there exists a subsequence
\((v_i)_{i \in M}\) of \((v_i)\) such that for every
\(x^* \in B_{X^*}\), the block basis
\(V_{x^*}= \{ v_i: i \in M, \, |x^*(v_i)| \geq 8 d \delta_{j_0}\}\),
is \(S_{p_{j_0}}\) admissible.

We next choose \(v_0\), a generic \((\delta_{j_0}, n_{j_0})\) average
of \((v_i)_{i \in M}\). It is easily seen that for some
\(x_0^* \in \mathcal{A}_{j_0}\) we have that
\(x_0^*(v_0) \geq \delta_{j_0}\). Therefore,
\(\|v_0\|_{j_0} \geq \delta_{j_0}\). On the other hand,
\(V_{x^*}\) is \(S_{p_{j_0}}\) admissible, for every
\(x^* \in B_{X^*}\) and \(p_{j_0} < n_{j_0}\).
It follows that \(\|v_0\| \leq (8d +1) \delta_{j_0}\).
We let \(v= \frac{v_0}{\|v_0\|}\) and observe that
\(\|v\|_{j_0} \geq \frac{1}{8d+1}\).

Let now \(j > j_0\). Arguing similarly, we can find 
a normalized block basis \((w_i)\) of \((u_i)\) and a 
generic \(({\delta_j}^2, n_j)\) average \(w_0\) of \((w_i)\) such
that \(v < w_0\) and \(\delta_j \leq \|w_0\| \leq 
(8d+1)\delta_j\). We let \(w=\frac{w_0}{\|w_0\|}\).
We are going to show that \(\|w\|_{j_0} \leq 
(8d+5) \delta_{j_0}\).
Suppose that \(\delta_{j_0} \sum_{i=1}^k x_i^*
\in \mathcal{A}_{j_0}\), and let \(E^*\) denote
the collection of the ranges of the \(x_i^*\)'s.
Let 
\(D= \{w_r: \, |\sum_{i=1}^k x_i^*(w_r)| \geq 8d \delta_j\}\).
Observe that by the choice of \((w_i)\) we have that
\(D(E^*,1)\) is \(2 S_{n_{j_0} + p_j}\) admissible.
On the other hand \(D(E^*,2)\) is \(2 S_{n_{j_0}}\)
admissible and thus \(D\) is \(4 S_{2p_j}\) admissible.
Because \(2p_j < n_j\), we obtain the estimate
\(\sum_{i=1}^k x_i^*(w_0) \leq (8d +4) \delta_j\).
Hence, \(\|w\|_{j_0} \leq (8d+5) \delta_{j_0}\), as claimed.
Finally, \(\frac{\|v\|_{j_0}}{\|w\|_{j_0}} \geq 
\frac{1}{(8d+1)(8d+5)\delta_{j_0}}\).
The proof is now complete since \(j_0\) was arbitrary.
\end{proof} 
\begin{Prop} \label{tin}
Let \(X_r\) have a shrinking basis \((e_k^r)_{k=1}^\infty\),
\(r=1,2\). Assume that
\(X_r\) satisfies the \((d_r,N_r,P_r,\mathbf{a})\)
distortion property, \(r=1,2\), and that
\(\mathbf{a}=(\delta_i)\) satisfies 
\(\lim_i \frac{\delta_{i+1}}{\delta_i} = 0\).
Suppose that for every \(i_0 \in \mathbb{N}\) there
exist \(i > j > i_0\) such that \(n_i^1 = n_j^2\),
where \(N_r=(n_k^r)_{k=1}^\infty\), \(r=1,2\).
Then \(X_1\) and \(X_2\) are totally incomparable.
\end{Prop}
\begin{proof}
Suppose the assertion is false. A standard perturbation
argument yields a normalized block basis \((u_k)\) of \((e_k^1)\)
equivalent to a block basis \((w_k)\) of \((e_k^2)\).
Let \(T\) be an isomorphism from \([(u_k)]\) onto
\([(w_k)]\) such that \(T(u_k)=w_k\), for all \(k \in \mathbb{N}\).
We can choose \(i_0 \in \mathbb{N}\) such that
\(\frac{\delta_{i+1}}{\delta_i} < \frac{1}{d_1 \|T\|\|T^{-1}\|}\),
for every \(i \geq i_0\).
Our assumptions allow us to choose \(i > j > i_0\) such that
\(n_i^1=n_j^2\). Let \((v_k)\) be a normalized block basis
of \((u_k)\) having no subsequence which is a
\(d_1 \delta_i\)-\(\ell_1^{n_i^1}\) spreading model.
But since \((T(v_k))\) is a block basis of \((w_k)\),
it follows that for every \(F \in S_{n_j^2}\) and all
choices of scalars \((a_k)_{k \in F}\)
\[\biggl \| \sum_{k \in F} a_k T(v_k) \biggr \| \geq
  \frac{\delta_j}{\|T^{-1}\|} \sum_{k \in F} |a_k|.\]
Hence, \(\| \sum_{k \in F} a_k v_k \| \geq 
\frac{\delta_j}{\|T\|\|T^{-1}\|} \sum_{k \in F} |a_k|\),
for every \(F \in S_{n_j^2}\) and all
choices of scalars \((a_k)_{k \in F}\).
However, \(\frac{\delta_i}{\delta_j} \leq 
\frac{\delta_{j+1}}{\delta_j}\), and therefore
\(\frac{\delta_j}{\|T\|\|T^{-1}\|} > d_1 \delta_i\).
Thus, \((v_k)\) is a \(d_1 \delta_i\)-\(\ell_1^{n_i^1}\) 
spreading model contrary to our assumptions.
\end{proof}
\begin{Def}
Let \(M=(m_i) \in [\mathbb{N}]\) such that
\(m_1 > 6\) and \(m_i^2 < m_{i+1}\), for all
\(i \in \mathbb{N}\). Choose \(L \in [\mathbb{N}]\),
\(L=(l_i)\) such that \(l_1 > 4\) and
\(2^{l_i} > m_i\), for all \(i \in \mathbb{N}\).
The infinite subset \(N=(n_i)\) of \(\mathbb{N}\)
is said to be \(M\)-good, if
\(l_j(f_j^N + 1) < n_j\), for all \(j \in \mathbb{N}\).
In the above, \((f_j^N)\) is the sequence given by
\(f_1^N=1\) while for \(j \geq 2\),
\[ f_j^N = \max \biggl \{ \sum_{i < j} \rho_i n_i :
\, \rho_i \in \mathbb{N} \cup \{0\} \, (i < j), \,
   \prod_{i < j} m_i^{\rho_i} < m_j^3 \biggr \}.\]
\end{Def}
Note that \(f_j^N\) is well defined because \(m_1 > 1\).
It is easy to see that for every \(P \in [\mathbb{N}]\)
there exists \(N \in [P]\) which is \(M\)-good.
The main result of Section \ref{S:5} is the following
\begin{Thm} \label{dhi}
Suppose \(N=(n_i)\) is \(M\)-good. Set 
\(N^{(2)}=(n_{2i})\), \(F^{(2)}=(f_{2i}^N + 2)\)
and \(\mathbf{a}= (\frac{1}{m_{2i}})\). Then
there exists a reflexive H.I. space \(X(N)\)
satisfying the \((6,N^{(2)},F^{(2)},\mathbf{a})\)
distortion property.
\end{Thm}
The proof is given in Section \ref{S:5}.
We now pass to the
\begin{proof}[{\bf Proof of Theorem \ref{T:1}}]
We first choose \(N_0 \in [\mathbb{N}]\) such that
every \(N \in [N_0]\) is \(M\)-good. To see that
such a \(N_0\) exists, set
\[\mathcal{D}= \{N \in [\mathbb{N}]: \, N \text{ is }
                   M-\text{good }\}.\]
We can easily verify that \(\mathcal{D}\) is closed
in the topology of pointwise convergence in \([\mathbb{N}]\),
and therefore it is a Ramsey set. Because 
\(\mathcal{D} \cap [R] \ne \emptyset\), for every
\(R \in [\mathbb{N}]\), the infinite Ramsey theorem
yields \(N_0 \in [\mathbb{N}]\) such that
\([N_0] \subset \mathcal{D}\), as claimed.

It is a well known fact that \([N_0]\) endowed 
with the topology of pointwise convergence is a perfect
Polish space. We let \([N_0]^2= [N_0] \times [N_0]\)
and set
\[G= \{(N,R) \in [N_0]^2, \, N=(n_i), \, R=(r_i)| \,
       \forall i_0 \in \mathbb{N}, \, \exists i > j > i_0 :
       \, n_{2i}=r_{2j}\}.\]
A straightforward application of the Baire category theorem
yields that \(G\) is a dense \(G_\delta\) subset of
\([N_0] \times [N_0]\). By a result of Kuratowski \cite{ku}
and Mycielski \cite{my} (cf. \cite{ke}, p. 129, Theorem
19.1, or Proposition 3.6 of \cite{gl}), there exists
\(C \subset [N_0]\) homeomorphic to the Cantor set such that
\((N_1, N_2) \in G\), whenever \(N_1\), \(N_2\) are distinct
elements of \(C\).

We can now apply Theorem \ref{dhi} to obtain
a family \(\{X(N): \, N \in C\}\) of reflexive H.I.
spaces such that for every \(N \in C\), \(X(N)\)
satisfies the \((6,N^{(2)},F^{(2)},\mathbf{a})\)
distortion property, where \(N^{(2)}\), \(F^{(2)}\)
and \(\mathbf{a}\) are as in the statement of Theorem
\ref{dhi}. Since \((N_1, N_2) \in G\) whenever
\(N_1\) and \(N_2\) are distinct elements of \(C\),
Proposition \ref{tin} implies that
\(X(N_1)\) and \(X(N_2)\) are totally incomparable.
The proof of the theorem is now complete.
\end{proof}
To construct H.I. spaces we shall make use of
the following
\begin{Thm} \label{hi}
Let \(X\) be a Banach space with a basis \((x_i)\).
Let \((n_j)\), \((k_j)\) be increasing sequences of
positive integers such that \(k_j < n_j\), for all
\(j \in \mathbb{N}\), and let \((\delta_j)\) be
a null sequence of positive scalars. Assume that
for every block subspace \(Y\) of \(X\) and every
\(j \in \mathbb{N}\) there exists a 
block basis \(z_1 < \cdots < z_p\) of \((x_i)\) in \(Y\)
such that letting \(t_i=\min \mathrm{supp} \, z_i\),
\(i \leq p\), the following are satisfied:
\begin{enumerate}
\item \(\{t_i:\, i \leq p\}\) is a maximal \(S_{n_j}\)
set and \(\|\sum_{i=1}^p a_i z_i\| \geq c_1 \delta_j
\|\sum_{i=1}^p a_i e_{t_i}\|_{n_j}\), for every sequence
\((a_i)_{i=1}^p\) in \(\mathbb{R}^+\).
\item \(\|\sum_{i=1}^p a_i z_i\| \leq c_2
\|\sum_{i=1}^p a_i e_{t_i}\|_{Ck_j} + c_3 {\delta_j}^2 \),
for every sequence
\((a_i)_{i=1}^p\) in \(\mathbb{R}\) with 
\(\sum_{i=1}^p |a_i| \leq 1\),
\end{enumerate}
where \(c_1\), \(c_2\) and \(c_3\) are absolute positive constants. Then
\(X\) has no infinite unconditional sequence.
If moreover, given \(Y\), \(Z\) block subspaces of
\(X\) and \(j \in \mathbb{N}\), such a block basis
\((z_i)_{i=1}^p\) can be found with the additional
property that \(z_i \in Y\), if \(i\) is odd,
while \(z_i \in Z\), if \(i\) is even, then \(X\)
is H.I.
\end{Thm}
\begin{proof}
Let \((u_i)\) be an infinite block basis of \((x_i)\),
and let \(j \in \mathbb{N}\). 
Set \(P=\{p_i: \, i \in \mathbb{N}\}\), where 
\(p_i = \min \mathrm{supp} \, u_i\).
We can find \(R \in [P]\) such that
\(\|[n_j]_1^L\|_{k_j} < {\delta_j}^2 \), for every
\(L \in [R]\). Let \(Y=[u_i: \, p_i \in R]\).
Choose \(z_1 < \cdots < z_p\) in \(Y\), according
to the hypothesis. There exists \(L \in [R]\)
such that \(\{t_i:\, i \leq p\} = \mathrm{supp} \,
              [n_j]_1^L\).
Put \(a_i=[n_j]_1^L(t_i)\), \(i \leq p\), and note that
\((a_i)_{i=1}^p\) is non-increasing.
We now have that
\[\biggl \|\sum_{i=1}^p a_i  z_i \biggr \| \geq c_1 \delta_j
 \biggl \|\sum_{i=1}^p a_i e_{t_i}\biggr \|_{n_j}= c_1 \delta_j.\]
On the other hand,
\[\biggl \|\sum_{i=1}^p (-1)^i a_i  z_i \biggr \| \leq
  c_2 \biggl \|\sum_{i=1}^p a_i e_{t_i} \biggr \|_{k_j} + c_3
{\delta_j}^2,\]
as \((a_i)_{i=1}^p\) is non-increasing.
Hence,
\[\biggl \|\sum_{i=1}^p (-1)^i a_i  z_i \biggr \| \leq
   (c_2 + c_3) {\delta_j}^2 \leq \frac{c_2 +c_3}{c_1} \delta_j
 \|\sum_{i=1}^p a_i  z_i\|.\]
Since \(j\) was arbitrary, \((u_i)\) is not unconditional.
The moreover statement is immediate.
\end{proof}
\section{Mixed Tsirelson spaces} \label{S:4}
Recall that
if \(\mathcal{M}\) is a set of finitely supported
signed measures on \(\mathbb{N}\) which satisfies the following:
\begin{enumerate}
\item \(e_n^* \in \mathcal{M}\), for all \(n \in \mathbb{N}\),
where \(e_n^*\) denotes the point mass measure at \(n\).
\item \(\mathcal{M}\) is symmetric i.e., if \(\mu \in \mathcal{M}\)
then \(- \mu \in \mathcal{M}\),
\item \(\mathcal{M}\) is pointwise bounded, that is
\(\mu(\{n\}) \leq 1\), for every \(\mu \in \mathcal{M}\),
\item \(\mathcal{M}\) is closed under restriction to initial
segments i.e., if \(\mu \in \mathcal{M}\),
then \(\mu | \{1, \dots , n\} \in \mathcal{M}\),
\end{enumerate}
then one can define a norm \(\|\cdot\|_{\mathcal{M}}\)
on \(c_{00}\) in the following manner:
\[\|\sum_{i=1}^{\infty} a_i e_i \|_{\mathcal{M}} = \sup
\{ \sum_{i=1}^{\infty} a_i \mu (\{i\}): \mu \in \mathcal{M}\},\]
for every finitely supported scalar sequence \((a_i)\). Of course,
\((e_i)\) is the natural basis of \(c_{00}\).
Letting \(X_{\mathcal{M}}\) denote the completion of
\((c_{00}, \|\cdot\|_{\mathcal{M}})\), we see that \((e_n)\)
is a normalized, monotone basis for \(X_{\mathcal{M}}\).
In case \(\mu | J \in \mathcal{M}\), for every \(\mu \in 
\mathcal{M}\) and \(J \subset \mathbb{N}\), then \((e_n)\)
is \(1\)-unconditional and bimonotone.

The main result of this section is
\begin{Thm} \label{mix}
Suppose \(N\) is \(M\)-good. There exists
\(\mathcal{M}\), a set of finitely supported
signed measures on \(\mathbb{N}\) satisfying conditions
1-4, above, and such that the following properties are fulfilled:
\begin{enumerate}
\item \((e_n)\) is an \(1\)-unconditional,
shrinking, bimonotone basis for \(X_{\mathcal{M}}\).
\item \(X_{\mathcal{M}}\) satisfies the
\((6,N,P,\mathbf{a})\) distortion property, where
\(P=(f_i^N + 2)\) and \(\mathbf{a}=(\frac{1}{m_i})\).
\end{enumerate}
\end{Thm}
We first give the construction of \(\mathcal{M}\) and prove
a number of lemmas necessary for the proof of Theorem \ref{mix}.

{\bf Construction of \(\mathcal{M}\)}. 
Given \(M=(m_i)\), \(N=(n_i)\), with \(N\) being \(M\)-good,
we construct \(\mathcal{M}\), a set of signed measures on
\(\mathbb{N}\) in the following manner: Let
\begin{align}
\mathcal{D} = \bigl \{ (t_1, \dots , t_{3n}): n \in \mathbb{N}, \, 
  &t_{3i-2} \in M \, (i < n), \, t_{3n-2}=0, \notag \\
  &t_{3i-1} \in [\mathbb{N}]^{< \infty} \setminus
  \{\emptyset\}, \,
  t_{3i} \in \{-1,1\} \, (i \leq n) \bigr \}. \notag
\end{align}
Given \(F \in \mathcal{D}^{< \infty}\), \(F \ne \emptyset\),
we let \(\mathcal{T}_F\) denote the set of all tuples
of length divisible by \(3\) which are initial segments
of elements of \(F\). We can partially order the elements of
\(\mathcal{T}_F\) by initial segment inclusion and thus 
\(\mathcal{T}_F\) becomes a finite tree with terminal nodes
precisely the members of \(F\). Given \(\alpha \in \mathcal{T}_F\)
then \(m \in M\) is an {\em \(M\)-entry} of \(\alpha\), if \(m \in \alpha\). 
We shall denote the last three entries of \(\alpha\)
by \(m_\alpha\), \(I_\alpha\) and \(\epsilon_\alpha\) respectively.
A rooted tree \(\mathcal{T} = \mathcal{T}_F\) 
(a tree is rooted if it has a unique root),
is said to be {\em appropriate} provided the following
properties hold:
\begin{enumerate}
\item If \(\alpha \in \mathcal{T}\) is terminal, then
\(I_\alpha = \{p_\alpha\}\), for some \(p_\alpha \in \mathbb{N}\).
\item If \(\alpha  \in \mathcal{T}\) is non-terminal and 
\(m_\alpha = m_j\), for some \(j \in \mathbb{N}\), then
\((I_\beta)_{\beta \in D_\alpha}\) is \(S_{n_j}\)-admissible
and \(I_\alpha = \cup_{\beta \in D_\alpha} I_\beta\). 
Here \(D_\alpha\) stands for the set of the immediate successors
of \(\alpha\) in \(\mathcal{T}\).
\end{enumerate}
We set
\[\mathcal{G} = \{ \mathcal{T}: \mathcal{T} \text{ is an
appropriate tree }\}.\]
We make the convention that the empty tree belongs to
\(\mathcal{G}\).
\begin{notation}
Let \(\mathcal{T} \in \mathcal{G}\) and \(\alpha \in \mathcal{T}\).
\begin{enumerate}
\item \(\alpha^{-}\) stands for the predecessor of \(\alpha\)
in \(\mathcal{T}\). In case \(\alpha\) is the root of 
\(\mathcal{T}\) we put \(\alpha^{-} = \emptyset\).
\item \(|\alpha|\) is the length of \(\alpha\). Thus,
\(|\alpha|=3n\) if \(\alpha = (t_1, \dots , t_{3n})\).
We now define \(o(\mathcal{T})= \max \{ |\beta|: \beta \in \mathcal{T}\}\),
the height of the tree \(\mathcal{T}\).
\item \(m(\alpha)= \prod_{m_i \in \alpha^{-}} m_i\). We 
set \(m(\alpha)=1\) if \(|\alpha|=3\).
\item \(n(\alpha) = \sum_{m_i \in \alpha^{-}} n_i\). We
set \(n(\alpha)=0\), if \(|\alpha|=3\).
\end{enumerate}
\end{notation}
Given \(\mathcal{T} \in \mathcal{G}\),   
set \[\mu_{\mathcal{T}}= \sum_{\alpha \in max \mathcal{T}}
      m(\alpha)^{-1} \epsilon(\alpha) \epsilon_\alpha 
      e_{p_\alpha}^*,\]
where \(max \mathcal{T}\) is the set of terminal nodes of
\(\mathcal{T}\) and \(I_\alpha = \{p_\alpha\}\) for
\(\alpha \in max \mathcal{T}\). We have also set
\(\epsilon(\alpha)= \prod_{\beta < \alpha}
\epsilon_\alpha \) for \(\alpha \in \mathcal{T}\).   
We make the convention \(\epsilon(\alpha)=1\), if \(|\alpha|=3\).
We also set \(\mu_{\emptyset}=0\).
Of course, \(\mu_\mathcal{T}\) is a finitely supported signed
measure on \(\mathbb{N}\) whose support is equal
to \(I_{\alpha_0}\), where \(\alpha_0\) is the root of
\(\mathcal{T}\). We also observe that \(|\mu_\mathcal{T}(\{n\})|
\leq 1\), for all \(n \in \mathbb{N}\).

We finally set \(\mathcal{M}= 
\{\mu_\mathcal{T}: \mathcal{T} \in \mathcal{G}\}\).
Note that \(e_n^* \in \mathcal{M}\) as
\(\bigl \{( 0, \{n\}, 1 ) \bigr \} \in \mathcal{G}\).
We shall introduce some more notation in order to investigate
properties of the set \(\mathcal{M}\).
\begin{notation}
Let \( \mathcal{T} \in \mathcal{G}\) and let \(\alpha_0\)
denote its root.
\begin{enumerate}
\item Given \(\alpha \in \mathcal{T}\) set
\(\mathcal{T}_\alpha = \{\beta \setminus \alpha^{-}: \beta \in
\mathcal{T}, \alpha \leq \beta\}\). 
Clearly, \(\mathcal{T}_\alpha \in \mathcal{G}\).
\item We let \(w(\mathcal{T})=1\), if \(|\mathcal{T}|=1\).
In case \(m_{\alpha_0} \in M\), we set
\(w(\mathcal{T})=m_{\alpha_0}\). 
\item Let \(J \subset \mathbb{N}\).
We let \(\mathcal{T}|J \) denote the tree resulting from
\(\mathcal{T}\) by keeping only those \(\alpha \in \mathcal{T}\) 
for which \(I_\alpha \cap J \ne \emptyset\) and replacing
\(I_\alpha\) by \(I_\alpha \cap J\).
It is easy to see that 
\(\mathcal{T}|J \in \mathcal{G}\).
\item We let \(- \mathcal{T}\) 
denote the tree resulting from \(\mathcal{T}\) by
changing \(\epsilon_{\alpha_0}\) to \(-\epsilon_{\alpha_0}\).
Clearly, \(- \mathcal{T} \in \mathcal{G}\) and moreover
\(\mu_{- \mathcal{T}}= - \mu_\mathcal{T}\).
\end{enumerate}
\end{notation} 
\begin{remark}
Let \(\mathcal{T} \in \mathcal{G}\).
\begin{enumerate}
\item If \(J \subset \mathbb{N}\), then
\(\mu_{\mathcal{T} | J} = \mu_\mathcal{T} | J\).
\item If \(\alpha \in \mathcal{T}\) then
\(m(\alpha) \epsilon(\alpha) \mu_\mathcal{T} | I_\alpha =
\mu_{\mathcal{T}_\alpha}\).
\end{enumerate}
\end{remark}
\begin{remark}
Suppose \(\mathcal{T}_i \in \mathcal{G}\), \(i \leq n\).
Let \(\alpha_i\) be the root of \(\mathcal{T}_i\), \(i \leq n\).
We shall say that \(\{\mathcal{T}_i: i \leq n\}\) is \(S_\xi\)-admissible,
\(\xi < \omega\), if \(\{I_{\alpha_i}: i \leq n\}\) is.
We shall also write \(\mathcal{T}_1 < \cdots < \mathcal{T}_n\)
if \(I_{\alpha_1} < \cdots < I_{\alpha_n}\).
It is easy to see that if \(\mathcal{T}_1 < \cdots < \mathcal{T}_n\)
is \(S_{n_j}\)-admissible then 
\(\frac{\sum_{i=1}^n \mu_{\mathcal{T}_i}}{m_j} \in \mathcal{M}\).
\end{remark}
It follows by our preceding remarks
that \(\mathcal{M}\) is pointwise bounded, symmetric
and closed under restriction to subsets of \(\mathbb{N}\).
Hence \((e_n)\) is an \(1\)-unconditional, bimonotone basis for
\(X_{\mathcal{M}}\). 
It is not hard to check that \(X_{\mathcal{M}}\) is
isometric to 
\(T(\frac{1}{m_i}, S_{n_i})_{i=1}^{\infty}\).
We also obtain by our preceding remarks
that if \((x_i)_{i=1}^k\)
is an \(S_{n_j}\)-admissible block basis of \((e_n)\)
then \(\|\sum_{i=1}^k x_i\| \geq \frac{1}{m_j} \sum_{i=1}^k
\|x_i\|\). Hence \(X_{\mathcal{M}}\) is an asymptotic 
\(\frac{1}{m_j}\)-\(\ell_1^{n_j}\) space.
It follows that \((e_n)\) is boundedly complete.
Let now \(\nu\)
be a \(w^*\)-cluster point of \(\mathcal{M}\).
Using the reflexivity argument of \cite{ad} ( cf. also \cite{ts}),
one obtains that for every \(\epsilon > 0\) there exists
\(k \in \mathbb{N}\) such that 
\(\|\nu | \, [e_i : \, i \geq k] \| < \epsilon\).
It follows from this that \((e_n)\) is shrinking and thus
\(X_{\mathcal{M}}\) is reflexive.
\begin{remark}
Suppose \((u_n)\) is a normalized block basis of
\((e_n)\) and \(u\) an \((\epsilon,n_j)\) average
of \((u_n)\). Then \(\frac{1}{m_j} \leq \|u\| \leq 1\).
\end{remark}
\begin{Lem} \label{L:1}
Let \(\mathcal{T} \in \mathcal{G}\).
Let \(F\) be a subset of
\(\mathcal{T}\) consisting of pairwise incomparable nodes.
Then \(\{I_\alpha: \alpha \in F\}\)
is \(S_p\)-admissible, where \(p = \max \{n(\alpha): \alpha \in F\}\).
\end{Lem}
\begin{proof}
By induction on \(o(\mathcal{T})\). If \(o(\mathcal{T})=3\)
the assertion of the lemma is trivial. Assuming the assertion
true when \(o(\mathcal{T}) < 3k\), \(k > 1\), let \(\mathcal{T}
\in \mathcal{G}\) with \(o(\mathcal{T})=3k\).  
If \(|F|=1\) there is nothing to prove. So assume \(|F|
\geq 2\). Let \(\alpha_0\) be the root of \(\mathcal{T}\)
and let \(w(\mathcal{T})=m_i\) for some \(i \in \mathbb{N}\).
We denote by \(D\)
the set of immediate successors of \(\alpha_0\) in
\(\mathcal{T}\). 
Given \(\alpha \in D\) let
\(F_\alpha = \{\beta \in F: \alpha \leq \beta\}\).
Because \( o(\mathcal{T}_\alpha) \leq 3k-3\) 
we can apply the induction hypothesis
on \(\mathcal{T}_\alpha\) and the set \(\{\beta \setminus \alpha^{-}:
\beta \in F_\alpha\}\) to deduce that the collection
\(\{I_\beta: \beta \in F_\alpha\}\) is \(S_{p_1}\)-admissible,
where \(p_1= \max \{ n(\beta \setminus \alpha^{-}):
\beta \in F_\alpha\}\). Since
\(n(\beta \setminus \alpha^{-})=n(\beta)-n(\alpha)\)
and \(n(\alpha)=n_i\) whenever \(\alpha \in D\),
we obtain that 
\(\{I_\beta: \beta \in F_\alpha\}\) is
\(S_{p-n_i}\)-admissible, for every \(\alpha \in D\).
But also, \(\{I_\alpha: \alpha \in D\}\) is \(S_{n_i}\)-admissible
whence \(\{I_\alpha: \alpha \in F\}\) is \(S_p\)-admissible.
\end{proof}
To simplify our notation, we set \(f_j=f_j^N\).
We make the following observation: Let \(\mathcal{T} \in \mathcal{G}\)
and let \(\alpha \in \mathcal{T}\). Assume that \(m(\alpha) < m_j^3\)
and that all \(M\)-entries of \(\alpha^{-}\) are smaller than 
\(m_j\). Then \(n(\alpha) \leq f_j\).
Our next lemma will be crucial for the proof of the main result.
\begin{Lem}[Decomposition Lemma] \label{L:2}
Let   
\(\mathcal{T}_0 \in \mathcal{G}\).
Let \(j \in \mathbb{N}\) such that 
\(w(\mathcal{T}_0) < m_j\). Then 
there exist an \(S_{f_j}\)-admissible subset \(\mathcal{G}_0\)
of \(\mathcal{G}\) and a scalar sequence 
\((\lambda_\mathcal{T})_{\mathcal{T} \in \mathcal{G}_0}\) in \([-1,1]\)
so that the following are satisfied:
\begin{enumerate}
\item \(\mu_{\mathcal{T}_0}= 
      \sum_{\mathcal{T} \in \mathcal{G}_0} \lambda_\mathcal{T} 
      \mu_\mathcal{T}\).
\item For each \(\mathcal{T} \in \mathcal{G}_0\) at least one of 
the following hold: either \(w(\mathcal{T})=1\) (thus
\(\mu_\mathcal{T} = \pm e_{\mathcal{T}(p)}^*\) for some
\(\mathcal{T}(p) \in \mathbb{N}\)), or 
\(w(\mathcal{T}) \geq m_j\), or
\(|\lambda_\mathcal{T}| \leq \frac{1}{m_j^2}\).
\end{enumerate}
\end{Lem}  
\begin{proof}
Let \(\mathfrak{B}\) denote the set of all branches
of \(\mathcal{T}_0\) (a branch is a maximal well ordered subset
of \(\mathcal{T}_0\)).
If \(w(\mathcal{T}_0)=1\) the assertion is trivial. So assume
that \(w(\mathcal{T}_0)=m_{i_0}\) for some \(i_0 < j\).
Given \(b \in \mathfrak{B}\) set
\[\alpha^1(b)= \max \{ \beta \in b: m(\beta) < m_j^2
\text{ and if } m_i \in \beta^{-} \text{ then } i < j\}.\]
Note that \(\alpha^1(b)\) is well defined 
and that \((m_{i_0}, I, \epsilon) < \alpha^1(b)\) since \(i_0 < j\)
(\((m_{i_0},I, \epsilon)\) being the root of \(\mathcal{T}_0\)).

Let us say that \(b \in \mathfrak{B}\) is of type \(1\) if
\(\alpha^1(b)\) is terminal in \(\mathcal{T}_0\).
If \(b\) is not of type \(1\) then it is of type
\(2\) (resp. \(3\)), if the last \(M\)-entry of \(\alpha^1(b)\)
is greater than or equal (resp. smaller than) \(m_j\).
We then denote by \(\alpha^2(b)\) the immediate successor
of \(\alpha^1(b)\) in \(b\).

We let \(A_1 = \{\alpha^1(b): b \in \mathfrak{B}
\text{ is of type } 1\}\),
\(A_2 = \{\alpha^1(b): b \in \mathfrak{B}
\text{ of type } 2 \}\) and 
\(A_3 = \{\alpha^2(b): b \in \mathfrak{B}
\text{ is of type } 3\}\).
Observe that the following properties hold:
\begin{enumerate}
\item If \(\alpha \in A_3\) then all \(M\)-entries of \(\alpha^{-}\)
are smaller than \(m_j\), \(m(\alpha^{-}) < m_j^2\), yet
\( m_j^2 \leq m(\alpha) < m_j^3\).
\item If \(\alpha \in A_2\), then \(\alpha\) is non-terminal,
all \(M\)-entries in \(\alpha^{-}\) are smaller than \(m_j\),
the last \(M\)-entry of \(\alpha\) is greater than or equal to
\(m_j\) and \(m(\alpha) < m_j^2\).
\item If \(\alpha \in A_1\) then \(\alpha\) is terminal,
all \(M\)-entries in \(\alpha^{-}\) are smaller than \(m_j\)
and \(m(\alpha) < m_j^2\).
\end{enumerate}
It is not hard to check now that \(A= \cup_{t=1}^3 A_t \)
consists of pairwise incomparable nodes of 
\(\mathcal{T}_0\) and hence 
\(\{I_\alpha: \alpha \in A\}\)
consists of successive subsets of \(\mathbb{N}\).
Moreover, \(I = \cup \{
I_\alpha: \alpha \in A\}\).
Because \(m(\alpha) < m_j^3\) and all \(M\)-entries of 
\(\alpha^{-}\) are smaller than \(m_j\) whenever
\(\alpha \in A\), we obtain that
\(n(\alpha) \leq f_j\) for all \(\alpha \in A\).
Lemma \ref{L:1} now yields that
\(\{I_\alpha: \alpha \in A\}\)
is \(S_{f_j}\)-admissible.
Finally, we let \(\mathcal{G}_0= \{ (\mathcal{T}_0)_{\alpha}:
\alpha \in A\}\).  
Since \(m(\alpha) \epsilon(\alpha) \mu_{\mathcal{T}_0} | I_\alpha =
\mu_{(\mathcal{T}_0)_\alpha}\), for all \(\alpha \in \mathcal{T}\),
we set \(\lambda_{(\mathcal{T}_0)_\alpha}=
\frac{1}{m(\alpha) \epsilon(\alpha)}\) for
\(\alpha \in A\).
We can easily verify that the desired properties hold.
\end{proof} 
In the sequel, we shall be using a variety of block bases of
\((e_n)\). The support of each of them will always be taken
with respect to \((e_n)\).
\begin{Lem} \label{L:4}
Let \((u_n)\) be a normalized block basis of
\((e_n)\).
Let \(j \in \mathbb{N}\), \(j \geq 2\) and let \(u\) be a generic 
\((\epsilon, f_j + 1)\) average of \((u_n)\) with
\(\epsilon < \frac{1}{2 m_j}\). Let \(i < j\)
and let \(\mathcal{T}_1 < \cdots < \mathcal{T}_t\) in \(\mathcal{G}\)
be \(S_{n_i}\)-admissible.
Then \(\sum_{k=1}^t \mu_{\mathcal{T}_k} (u) \leq 2\).
In particular, \(\mu_\mathcal{T} (u) \leq \frac{2}{w(\mathcal{T})}\),
if \(w(\mathcal{T}) < m_j\).
\end{Lem}
\begin{proof}
Observe that \(\frac{1}{m_j} \sum_{k=1}^t \mu_{\mathcal{T}_k}
\in \mathcal{M}\) and hence 
\(\sum_{k=1}^t \mu_{\mathcal{T}_k} (u_n) \leq m_j\),
for all \(n \in \mathbb{N}\).
Let \(P=(p_n)\), where \(p_n = \min \mathrm{supp} \, u_n\).
Set \(\xi = f_j + 1\) and 
suppose that \(u= \sum_{n=1}^\infty \xi_1^R(p_n) u_n\),
for some \(R \in [P]\). Let \(E^*\) denote the collection
of the ranges of the \(\mu_{\mathcal{T}_k}\)'s, and let \(D\)
denote the collection of those \(u_n\)'s whose support
intersects at least one member of 
\(E^*\). Put \(I_r=\{n \in \mathbb{N}:
u_n \in D(E^*,r)\}\), \(r=1,2\). Because \(D(E^*,2)\)
is \(2S_{n_i}\)-admissible and \(n_i \leq f_j\), we obtain that
\(\sum_{k=1}^t \mu_{\mathcal{T}_k}(\sum_{n \in I_2} \xi_1^R(p_n) u_n)
\leq m_j 2 \epsilon\). On the other hand we clearly have that
\(\sum_{k=1}^t \mu_{\mathcal{T}_k}(\sum_{n \in I_1} \xi_1^R(p_n) u_n)  
\leq 1\). Thus, \(\sum_{k=1}^t \mu_{\mathcal{T}_k}(u) \leq 2\). 
\end{proof}
\begin{Lem} \label{L:5}
Let \((u_n)\) be a normalized block basis of
\((e_n)\). Let \(\epsilon > 0\) and \(j \in \mathbb{N}\).
Then there exists a smoothly normalized 
\((\epsilon, f_j + 1)\) average of \((u_n)\).
\end{Lem}
\begin{proof}
Let \(P=(p_n)\), where \(p_n= \min \mathrm{supp} \, u_n\) for
\(n \in \mathbb{N}\). We can assume without loss of generality
that \(\|\xi_1^R\|_{\xi - 1} < \epsilon\) for every
\(R \in [P]\) where \(\xi=f_j + 1\).
We are going to show that there exists a normalized block basis
of \((u_n)\) admitting a generic
\((\epsilon,\xi)\) average of norm at least \(\frac{1}{2}\).
Suppose instead that this were false. Then it is easy to
construct for every \(1 \leq r \leq l_j \), a block basis
\((u_i^r)\) of \((u_i)\) so that letting \(p_i^r =
\min \mathrm{supp} \, u_i^r\) and \(P_r = ( p_i^r )\) the following
are satisfied:
\begin{enumerate}
\item \((u_i^r)\) is a block basis of \((u_i^{r-1})\).
(\(u_i^0=u_i\))
\item \(u_i^r = \sum_{n=1}^\infty \xi_i^{P_{r-1}} 
(p_n^{r-1}) \frac{u_n^{r-1}}{\|u_n^{r-1}\|}\), for all
\(i \in \mathbb{N}\).
(\(p_n^0=p_n\))
\item \(\|u_i^r\| < \frac{1}{2}\), for all \(i \in \mathbb{N}\).
\item For every \(i \in \mathbb{N}\),
if \(u_i^r = \sum_{n \in F_i^r} a_n u_n\) with \(a_n > 0\)
for \(n \in F_i^r\),
then \(\sum_{n \in F_i^r} a_n \geq  2^{r-1}\) and 
\((u_n)_{n \in F_i^r}\) is \(S_{\xi r}\)-admissible.
\end{enumerate}
The construction is easily done by induction.
Taking \(r=l_j \) we see from 3. that
\(\|u_i^{l_j}\| < \frac{1}{2}\). On the other
hand 4. implies that \(\|u_i^{l_j}\| \geq \frac{2^{l_j -1}}{m_j}\)
as \(\xi l_j < n_j\). Thus, \(m_j > 2^{l_j}\) contradicting
the choice of \(l_j\).
\end{proof}
Our next lemma yields that \(X_{\mathcal{M}}\) satisfies the
\((6,N,F,\mathbf{a})\) distortion property where
\(F=(f_i+2)\) and \(\mathbf{a}=(\frac{1}{m_i})\).
\begin{Lem} \label{L:6}
Let \((u_j)\) be a normalized block basis of \((e_j)\).
Suppose that \((y_j)\)
is a block basis of \((u_j)\)
so that \(y_j\) is a smoothly normalized \((\epsilon_j , f_j + 1)\)
average of \((u_j)\) with \(\epsilon_j < \frac{1}{2m_j}\).
Given \(j_0 \in \mathbb{N}\) and \(J_0 \in [\mathbb{N}]\),
there exists \(J \in [J_0]\) such that \(j_0 < \min J\)
and for every \(\mathcal{T} \in \mathcal{G}\),
\(D_\mathcal{T} = \{y_j : \, j \in J, |\mu_\mathcal{T} (y_j)|
\geq \frac{5}{m_{j_0}} \}\) is
\(S_{f_{j_0} + 1}\)-admissible. 
\end{Lem}
\begin{proof}
Note first that Lemma \ref{L:5} guarantees the existence
of the block basis \((y_j)\).
Let \(P= (p_j)_{j \in J_0}\), where \(p_j = \min \mathrm{supp} \, y_j\).
By passing to a subsequence of \((y_j)_{j \in J_0}\), if necessary,
we can assume that the union of any \(4\) \(S_{f_{j_0}}\)
subsets of \(P\) belongs to \(S_{f_{j_0} + 1}\).
Choose \(J \in [J_0]\), \(J=(j_i)\), such that
\(j_0 < j_1\) and 
\(\|y_{j_i}\|_{\ell_1} < \frac{m_{j_{i+1}}}{m_{j_i}}\),
for every \(i \geq 2\) (if \(v=\sum_{i=1}^n a_i e_i\), then
\(\|v\|_{\ell_1}= \sum_{i=1}^n |a_i|\)).

Let \(\mathcal{T}_0 \in \mathcal{G}\). Suppose first that
\(w(\mathcal{T}_0) \geq m_{j_0}\). We show that in this case 
\(|D_{\mathcal{T}_0}| \leq 1\).
Indeed, suppose first that \(w(\mathcal{T}_0) < m_{j_1}\). 
Lemma \ref{L:4} yields that \(|\mu_{\mathcal{T}_0} (y_j)|
\leq \frac{4}{m_{j_0}}\) for all \(j \in J\)
whence \(D_{\mathcal{T}_0} = \emptyset\).

If \(w(\mathcal{T}_0) \geq m_{j_1}\) choose 
\(s \geq 2\) so that
\(m_{j_{s-1}} \leq w(\mathcal{T}_0) < m_{j_s}\).
Observe that if \(1 \leq i < s-1\) then
\[ |\mu_{\mathcal{T}_0} (y_{j_i})| \leq \frac{1}{w(\mathcal{T}_0)}
\|y_{j_i}\|_{\ell_1} 
< \frac{1}{w(\mathcal{T}_0)} \frac{m_{j_{s-1}}}{m_{j_i}}
< \frac{1}{m_{j_0}}.\]
When \(i \geq s\), Lemma \ref{L:4} yields
\(| \mu_{\mathcal{T}_0} (y_{j_i})| \leq \frac{4}{w(\mathcal{T}_0)}\)
\(< \frac{4}{m_{j_0}}\). Hence 
\(D_{\mathcal{T}_0} \subset \{y_{j_{s-1}}\}\) and so our claim holds.

The final case to consider is that of \(w(\mathcal{T}_0) < m_{j_0}\).
Clearly, \(D_{\mathcal{T}_0}=\emptyset\), if \(w(\mathcal{T}_0)=1\).
We employ the decomposition Lemma \ref{L:2} to find an
\(S_{f_{j_0}}\) admissible subset \(\mathcal{G}_0\) of \(\mathcal{G}\)
and scalars \((\lambda_\mathcal{T})_{\mathcal{T} \in \mathcal{G}_0}\)
satisfying the conclusion of Lemma \ref{L:2}.
Let \(E^*\) denote the collection
of the ranges of the \(\mu_{\mathcal{T}}\)'s (\(\mathcal{T} \in \mathcal{G}_0\)).
Our previous work implies that \(D_{\mathcal{T}_0}(E^*,1)\) is
\(2S_{f_{j_0}}\) admissible. But also, \(D_{\mathcal{T}_0}(E^*,2)\)
is \(2S_{f_{j_0}}\) admissible since \(\mathcal{G}_0\) is 
\(S_{f_{j_0}}\) admissible. It follows that
\(D_{\mathcal{T}_0}\) is \(S_{f_{j_0} + 1}\) admissible.
\end{proof}
\begin{proof}[{\bf Proof of Theorem \ref{mix}}]
Let \(U\) be a block subspace of \(X_{\mathcal{M}}\)
spanned by the normalized block basis \((u_j)\) of \((e_j)\).
Let \(j_0 \in \mathbb{N}\) and choose a block basis
\((y_j)_{j \in J}\) of \((u_j)\) satisfying the conclusion
of Lemma \ref{L:6}. 
Applying Corollary 3.4 of \cite{ano} (cf. also Corollary 3.3
of \cite{g}), we obtain that for every subsequence
of \((y_j)_{j \in J}\) which is a 
\(\delta\)-\(\ell_1^{f_{j_0} +2}\) spreading model,
it must be the case that 
\(\delta \leq \frac{5}{m_{j_0}}\) and thus
\(\tau(U, \frac{6}{m_{j_0}}) < f_{j_0} +2\).
\end{proof}
{\bf Terminology.}
Let \(j_0\) and \((y_j)_{j \in J}\) satisfy the conclusion of
Lemma \ref{L:6}. Every normalized \((\epsilon, n_{j_0})\)
average \(u\) of \((u_j)_{j=1}^{\infty}\) of the form
\(u=\frac{v}{\|v\|}\), where \(v\) is a generic
\((\epsilon, n_{j_0})\) average of \((y_j)_{j \in J}\),
will be called a {\em normalized} \((\epsilon, n_{j_0})\)
{\em average of} \((u_j)_{j=1}^{\infty}\) {\em resulting
from Lemma \ref{L:6}}. Note that Lemmas \ref{L:5} and 
\ref{L:6} guarantee the existence of such averages for
every block basis \((u_j)_{j=1}^{\infty}\).
\begin{Cor} \label{L:7}
Let \((y_j)\) satisfy the assumptions of Lemma \ref{L:6}.
Given \(j_0 \in \mathbb{N}\) and \(J_0 \in [\mathbb{N}]\),
there exists \(J \in [J_0]\) such that \(j_0 < \min J\)
and for every \(\mathcal{T}_0 \in \mathcal{G}\), \(w(\mathcal{T}_0) \ne m_{j_0}\),
\(D_{\mathcal{T}_0} = \{y_j : \, j \in J, |\mu_\mathcal{T} (y_j)|
\geq \frac{5}{m_{j_0} m_e} \}\) is
\(S_{f_{j_0} + 1}\)-admissible, where we have set
\(m_e= \min \{m_{j_0}, w(\mathcal{T}_0)\}\).
\end{Cor}
\begin{proof}
We choose \(J_0\) and \(j_0\) as we did in the proof
of Lemma \ref{L:6}. Suppose first that \(w(\mathcal{T}_0) > m_{j_0}\). 
Because \(m_i^2 < m_{i+1}\), the argument in the proof of Lemma
\ref{L:6} shows that \(|D_{\mathcal{T}_0}| \leq 1\).

When \(w(\mathcal{T}_0) < m_{j_0}\), we apply the decomposition
Lemma \ref{L:2} to find an
\(S_{f_{j_0}}\) admissible subset \(\mathcal{G}_0\) of \(\mathcal{G}\)
and scalars \((\lambda_\mathcal{T})_{\mathcal{T} \in \mathcal{G}_0}\)
satisfying the conclusion of Lemma \ref{L:2}.
Note that if \(\mathcal{T} \in \mathcal{G}_0\) and \(w(\mathcal{T})=m_{j_0}\),
then \(|\lambda_\mathcal{T}| \leq \frac{1}{w(\mathcal{T}_0)}\) and thus for
all \(j \in J\),
\(|\lambda_\mathcal{T} \mu_\mathcal{T}(y_j)| \leq \frac{4}{m_{j_0} m_e}\),
by Lemma \ref{L:4}. Using the splitting argument of Lemma
\ref{L:6}, we conclude that \(D_{\mathcal{T}_0}\) is
\(S_{f_{j_0} + 1}\)-admissible.
\end{proof}
\begin{Cor} \label{L:8}
Let \(u\) be a normalized 
\((\epsilon, n_{j_0})\)
average of \((u_j)_{j=1}^{\infty}\) resulting
from Lemma \ref{L:6} with \(\epsilon \leq \frac{1}{12m_{j_0}^2}\).
Let \(\mathcal{G}_0\) be an \(S_{n_i}\)-admissible subset
of \(\mathcal{G}\), \(i < j_0\), such that \(m_{j_0} \notin \{w(\mathcal{T}):
\mathcal{T} \in \mathcal{G}_0\}\). Then,
\(|\sum_{\mathcal{T} \in \mathcal{G}_0} \mu_\mathcal{T} (u)| \leq
  \frac{6}{m_e}\), where
\(m_e = \min \{w(\mathcal{T}):
\mathcal{T} \in \mathcal{G}_0\} \cup \{m_{j_0}\}\).
\end{Cor}
\begin{proof}
Set \(\xi=n_{j_0}\). Let \(p_j =\min \mathrm{supp} \, y_j\),
\(j \in J\) and \(P=\{p_j: \, j \in J\}\). There exists
\(R \in [P]\) so that \(u=\frac{v}{\|v\|}\), 
where \(v= \sum_{j \in J} \xi_1^R(p_j) y_j\) and
\(\|\xi_1^R\|_{n_{j_0}-1} < \epsilon\).
Note that \(\|v\| \geq \frac{1}{m_{j_0}}\).
Applying a splitting argument similar to
that of Lemma \ref{L:6} and taking in account Corollary
\ref{L:7}, we obtain that
\(\{y_j: \, j \in J, 
|\sum_{\mathcal{T} \in \mathcal{G}_0} \mu_\mathcal{T} (y_j)| \geq
  \frac{5}{m_{j_0}m_e}\}\)
is \(3 S_{2f_{j_0} + 1}\)-admissible.
The assertion follows from Lemma \ref{L:4} and the fact
that \(2f_{j_0} + 1 < n_{j_0}\).
\end{proof}
\begin{remark}
It is easy to see that in case 
\(w(\mathcal{T})=1\), for all \(\mathcal{T} \in \mathcal{G}_0\),
one obtains the estimate 
\(|\sum_{\mathcal{T} \in \mathcal{G}_0} \mu_\mathcal{T} (u)| \leq 
 \frac{1}{m_{j_0}}\).
\end{remark}
\section{Hereditarily indecomposable spaces} \label{S:5}
This section is devoted to the proof of Theorem \ref{dhi}.
Recall that \(X\) is H.I. if and only if,
for every pair of subspaces \(Y\), \(Z\) 
of \(X\) and every \(\epsilon > 0\), there exist
\(y \in Y\) and \(z \in Z\) so that \(y \ne z\) and
\(\|y-z\| \leq \epsilon \|y + z\|\).

Let \(M \in [\mathbb{N}]\), \(M=(m_i)\) and let  
\(N \in [\mathbb{N}]\), \(N=(n_i)\), which is
\(M\)-good.
Let \(\mathcal{M}\) be the set of measures constructed in
the previous section by using the sets \(M\) and \(N\).
We shall choose \(\mathcal{N} \subset \mathcal{M}\)
so that the resulting space \(X_{\mathcal{N}}\) is a reflexive H.I.
space satisfying the conclusion of Theorem \ref{dhi}.

We can find an injection \[\sigma \colon
\{ \mathcal{T}_1 < \cdots < 
\mathcal{T}_n: n \in \mathbb{N}, \mathcal{T}_i 
\in \mathcal{G} \, (i \leq n)\}
\to \{m_{2j}: j \in \mathbb{N}\}\]
so that \(\sigma(\mathcal{T}_1 , \dots , \mathcal{T}_n) > 
w(\mathcal{T}_i)\),
for all \(i \leq n\).
\begin{Def}
\begin{enumerate}
\item An \(S_p\)-admissible sequence
\(\mathcal{T}_1 < \cdots < \mathcal{T}_n\) in \(\mathcal{G}\) is
said to be \(S_p\)-dependent, \(p \in \mathbb{N}\),
if \(w(\mathcal{T}_1)=m_{2j_1}\), for some \(j_1 > \frac{p}{2}\),
and \(\sigma(\mathcal{T}_1 , \dots , \mathcal{T}_{i-1})= w(\mathcal{T}_i)\),
for all \(2 \leq i \leq n\).
\item Let \(\mathcal{T}_1 < \cdots < \mathcal{T}_n\) in \(\mathcal{G}\),
\(p \in \mathbb{N}\) and \(\mathcal{G}_0 \subset \mathcal{G}\).
We shall say that \(\mathcal{T}_1 < \cdots < \mathcal{T}_n\)
admits an \(S_p\)-dependent extension in \(\mathcal{G}_0\),
if there exist \(l \in \mathbb{N}\),
\(k \in \mathbb{N} \cup \{0\}\) and 
an \(S_p\)-dependent sequence
\(\mathcal{R}_1 < \cdots < \mathcal{R}_{n+k}\)
in \(\mathcal{G}_0\) so that \(\mathcal{R}_{k+i}| 
[l, \infty ) = \mathcal{T}_i\), for all \(i \leq n\).
\item A subset \(\mathcal{G}_0\) of \(\mathcal{G}\) is said
to be self-dependent, if the following property is satisfied
for every \(\mathcal{T} \in \mathcal{G}_0\):
Let \(\alpha \in \mathcal{T}\) so that its last \(M\)-entry
equals \(m_{2j+1}\) for some \(j \in \mathbb{N}\).
Let \(D_\alpha\) denote the set of immediate successors
of \(\alpha\) in \(\mathcal{T}\). Then
\(\{\mathcal{T}_\beta: \beta \in D_\alpha\}\) admits an
\(S_{n_{2j+1}}\)-dependent extension in \(\mathcal{G}_0\).
\end{enumerate}
\end{Def}
\begin{Def}
We let \(\mathfrak{D}\) denote the union of all non-empty,
self-dependent, symmetric and closed under restriction
to intervals, subsets of \(\mathcal{G}\).
Recall that \(\mathcal{G}_0 \subset \mathcal{G}\) is symmetric if
\(- \mathcal{T} \in \mathcal{G}_0\) 
whenever \(\mathcal{T} \in \mathcal{G}_0\).
\(\mathcal{G}_0\) is closed under interval restrictions if
\(\mathcal{T} | J \in \mathcal{G}_0\) whenever
\(\mathcal{T} \in \mathcal{G}_0\) and \(J\) is an interval.
\end{Def}
Of course \(\mathfrak{D}\) is a maximal, under inclusion,
subset of \(\mathcal{G}\) with respect to the aforementioned
properties. Set \(\mathcal{N}= 
\{\mu_\mathcal{T}: \mathcal{T} \in \mathfrak{D}\}\).
We will show that \(X_{\mathcal{N}}\) is H.I.
\begin{remark}
The maximality of \(\mathfrak{D}\) implies the following:
\begin{enumerate}
\item \(e_n^* \in \mathcal{N}\), for all \(n \in \mathbb{N}\).
\item If \(\mathcal{T} \in \mathfrak{D}\),
then \(\mathcal{T}_\alpha \in \mathfrak{D}\), for all
\(\alpha \in \mathcal{T}\) and so the decomposition
Lemma \ref{L:2} holds for \(\mathfrak{D}\).
\item If 
\(\mathcal{T}_1 < \cdots < \mathcal{T}_k\) in \(\mathfrak{D}\)
is \(S_{n_{2i}}\)-admissible, \(i \in \mathbb{N}\), then
\(\frac{\mu_{\mathcal{T}_1} + \cdots + \mu_{\mathcal{T}_k}}{m_{2i}}
\in \mathcal{N}\).
\item If 
\(\mathcal{T}_1 < \cdots < \mathcal{T}_k\) in \(\mathfrak{D}\)
is \(S_{n_{2i+1}}\)-dependent, \(i \in \mathbb{N}\), then
\(\frac{\mu_{\mathcal{T}_1} + \cdots + \mu_{\mathcal{T}_k}}{m_{2i+1}}
\in \mathcal{N}\).
\item Because of 3., all the results obtained in
the previous section 
about \((\epsilon, \xi)\) averages in \(X_{\mathcal{M}}\),
where \(\xi\) is either \(n_j\) or \(f_j +1\) 
for some \(j \in \mathbb{N}\),
still hold in \(X_{\mathcal{N}}\) provided
\(j\) is even.  
\end{enumerate}
\end{remark}
Note that \(X_{\mathcal{N}}\) is reflexive by the same argument
that showed \(X_{\mathcal{M}}\) is reflexive.
Thus \((e_i)\) is a shrinking basis for \(X_{\mathcal{N}}\).
\begin{proof}[{\bf Proof of Theorem \ref{dhi}}]
It follows from Theorem \ref{mix} and our preceding remarks that
\(X_{\mathcal{N}}\) satisfies the 
\((6,N^{(2)},F^{(2)},\mathbf{a})\) distortion property.
We show that \(X_{\mathcal{N}}\) is H.I.
This is accomplished through Theorem \ref{hi}.
Let \((u_n)\) be a normalized block basis of \((e_n)\)
and let \(j \in \mathbb{N}\). Set
\(P=(p_n)\), where \(p_n= \min \, \mathrm{supp} \, u_n\).
We can assume that the union of any \(7\) 
\(S_{f_{2j + 1}}\) subsets of \(P\) belongs to
\(S_{f_{2j + 1} + 1}\).
Successive applications of Corollary \ref{L:8}
yield a normalized block basis \(g_1 < \cdots < g_p\)
of \((u_n)\), \(\mathcal{T}_1 < \cdots < \mathcal{T}_p\) in
\(\mathfrak{D}\), and integers \(j_1 < \cdots < j_p\)
with \(2j+1 < j_1\), satisfying the following:
\begin{enumerate}
\item \(g_i\) is a normalized \((\frac{1}{12m_{2j_i}^2}, n_{2j_i})\)
average of \((u_n)\) resulting from Lemma \ref{L:6}.
\item \(w(\mathcal{T}_i)= m_{2j_i}\), 
\(\mathrm{supp} \, \mu_{\mathcal{T}_i} \subset r(g_i)\) 
and 
\(\mu_{\mathcal{T}_i}(g_i) > \frac{1}{2}\), for all \(i \leq p\).
\item \(\sigma(\mathcal{T}_1, \dots , \mathcal{T}_{i-1})=w(\mathcal{T}_i)\),
for all \(i \leq p\).
\item \(\{g_i: \, i \leq p\}\) is maximally
\(S_{n_{2j + 1}}\)-admissible.
\end{enumerate}
Put \(\theta_i= (\mu_{\mathcal{T}_i}(g_i))^{-1}\), 
\(z_i=\theta_i g_i\),
and note that \(1 \leq \theta_i < 2\), \(i \leq p\).
We'll show that \((z_i)_{i=1}^p\) satisfies conditions
1. and 2. of Theorem \ref{hi}, with \(\delta_j= \frac{1}{m_{2j+1}}\),
``\(n_j \)''\(=n_{2j+1}\) and \(k_j=f_{2j + 1} + 1\).
Condition 1. is immediate
since \(\mathcal{T}_1 < \dots < \mathcal{T}_p\) 
is \(S_{n_{2j + 1}}\)-dependent.   
Condition 2. is achieved by establishing the following

{\bf Claim:} Given \(\mathcal{T} \in \mathfrak{D}\), there exist
intervals \(J_1 < \cdots < J_s\) in \(\{1, \dots, p\}\)
so that
\begin{enumerate}
\item \(\{z_{\min J_t}: \, t \leq s\}\) is 
\(S_{f_{2j + 1} + 1}\)-admissible.
\item \(\mu_\mathcal{T} \, | \, \{z_i: \, i \in J_t\}\) is constant
for all \(t \leq s\).
\item \(|\mu_\mathcal{T}(z_i)| < \frac{12}{m_{2j+1}^2}\), for all
\(i \notin \cup_{t=1}^s J_t\).
\end{enumerate}
To prove the claim suppose first that
\(w(\mathcal{T}) > m_{2j + 1}\).
Corollary \ref{L:8} yields that 
\(|\mu_\mathcal{T}(z_i)| \geq \frac{12}{m_{2j+1}^2}\),
for at most one \(i \leq p\), and thus the claim holds in
this case.

Next assume that \(w(\mathcal{T})=m_{2j+1}\).
Without loss of generality, 
there exist an \(S_{n_{2j + 1}}\)-dependent
sequence \(\mathcal{R}_1 < \cdots < \mathcal{R}_l\) in \(\mathfrak{D}\)
and an interval \(J\) so that 
\(\mu_\mathcal{T}= \frac{1}{m_{2j+1}} \sum_{k=1}^l \mu_{\mathcal{R}_k}
\, | \, J\).
Let \(i_0\) be the largest \(i\) for which
\(w(\mathcal{T}_i)\) is an element of
\(\{w(\mathcal{R}_k): \, k \leq l\}\), or let
\(i_0=0\), if no such \(i\) exists.
The injectivity of \(\sigma\) and Corollary 
\ref{L:8} imply that if \(i_0 \in \{0,1\}\),
or if \(w(\mathcal{T}_{i_0})=w(\mathcal{R}_1)\), then
\(|\mu_\mathcal{T}(z_i)| < \frac{12}{m_{2j+1}^2}\),
for all \(i \ne i_0\).

If \(i_0 > 1\), then the injectivity of \(\sigma\) yields
\(w(\mathcal{T}_{i_0})=w(\mathcal{R}_{i_0})\), 
\(\mathcal{T}_i = \mathcal{R}_i\)
for \(i < i_0\) yet \(\mathcal{T}_{i_0} \ne \mathcal{R}_{i_0}\).
It follows now by Corollary \ref{L:8}, that
\(|\mu_\mathcal{T}(z_i)| < \frac{12}{m_{2j+1}^2}\),
for all \(i > i_0\).
We also observe that there exists \(i_1 < i_0\)
such that \(\mu_\mathcal{T}(z_i)=0\), if \(i < i_1\),
while \(\mu_\mathcal{T}(z_i)=\frac{1}{m_{2j+1}}\)
if \(i_1 < i < i_0 - 1 \).
Concluding, there exist four intervals
\(J_1 < J_2 < J_3 < J_4\) in \(\{1, \dots , p\}\), 
some of which may possibly be empty, such that
\(\mu_\mathcal{T} \, | \, \{z_i: \, i \in J_t\}\)
is constant for every \(t \leq 4\), while
\(|\mu_\mathcal{T}(z_i)| < \frac{12}{m_{2j+1}^2}\),
for each \(i \notin \cup_{t=1}^4 J_t\).

Finally, assume \(w(\mathcal{T}) < m_{2j+1}\).
If \(w(\mathcal{T})=1\), the claim trivially holds
so suppose that \(w(\mathcal{T}) > 1\).
Choose \(\mathcal{G}_0 \subset \mathfrak{D}\) \(S_{f_{2j+1}}\)-admissible
and scalars \((\lambda_{\mathcal{R}})_{\mathcal{R} \in \mathcal{G}_0}\)
according to the decomposition Lemma \ref{L:2}.
By splitting the \(z_i\)'s into two sets, those whose support
intersects at least two of the ranges of the
\(\mu_{\mathcal{R}}\)'s, and those whose support intersects at most one,
we deduce from our previous work that there exist
intervals \(J_1 < \cdots < J_s\) in \(\{1,\dots,p\}\)
so that \(\{z_{\min J_t}: \, t \leq s\}\) is 
\(7 S_{f_{2j + 1}}\)-admissible,
\(\mu_\mathcal{T} \, | \, \{z_i: \, i \in J_t\}\) is constant
for all \(t \leq s\), and 
\(|\mu_\mathcal{T}(z_i)| < \frac{12}{m_{2j+1}^2}\), for all
\(i \notin \cup_{t=1}^s J_t\).
Thus the claim holds and the proof is complete.
\end{proof}

\end{document}